\documentclass[a4paper,11pt]{amsart}
\usepackage{amssymb}
\newtheorem{theorem}{Theorem}[section]
\newtheorem{lemma}[theorem]{Lemma}
\newtheorem{prop}[theorem]{Proposition}
\newtheorem{corollary}[theorem]{Corollary}

\theoremstyle{definition}
\newtheorem{definition}[theorem]{Definition}
\newtheorem{example}[theorem]{Example}

\theoremstyle{remark}
\newtheorem{remark}[theorem]{Remark}
\numberwithin{equation}{section}
\newcommand{\BB}{{\mathbb B}}
\newcommand{\DD}{{\mathbb D}}
\newcommand{\OO}{{\mathcal O}}

\newcommand{\calm}{{\mathcal M}}
\newcommand{\NN}{{\mathbb N}}

\newcommand{\CC}{{\mathbb C}}

\newcommand{\eps}{\varepsilon}

\DeclareMathOperator{\psh}{PSH}

\DeclareMathOperator{\dist}{dist}

\renewcommand{\phi}{\varphi}

\begin{document}
\title[Determination of the pluripolar hull]
{Determination of the pluripolar hull of graphs of certain holomorphic functions}

\author{Armen Edigarian}

\address{Institute of Mathematics, Jagiellonian University,
Reymonta 4/526, 30-059 Krak\'ow, Poland}
\email{edigaria@im.uj.edu.pl}
\thanks{The first author was supported in part by the KBN grant
No. 5 P03A 033 21. The first author is a fellow of the
Rector's Scholarship Fund at Jagiellonian University}

\author{Jan Wiegerinck}

\address{Faculty of Mathematics, University of Amsterdam,
Plantage Muidergracht 24, 1018 TV, Amsterdam, The Netherlands}
\email{janwieg@science.uva.nl}

\date{16 July 2003}




\begin{abstract}Let $A$ be a closed polar subset of a domain $D$ in $\CC$.
We give a complete description of the pluripolar hull
$\Gamma^*_{D\times\CC}$ of the graph $\Gamma$ of a holomorphic function definedon $D\setminus A$. To achieve this, we prove for pluriharmonic measure certain semi-continuity properties and a localization principle.
\end{abstract}

\maketitle

\section{Introduction}
Let $f$ be a holomorphic function on its {\em domain of existence} $D\subset\CC$ and
let
$\Gamma_f$ be its graph in $D\times \CC$. Answering a question of Levenberg, 
Martin and Poletsky \cite{LMP}, we showed in \cite{EW1} that it is possible that 
$\Gamma_f$ is {\em not} a complete pluripolar subset of $\CC^2$, but that the pluripolar 
hull of $\Gamma_f$ is strictly larger than $\Gamma_f$. In a subsequent paper \cite{EW2} 
we studied the pluripolar hull $(\Gamma_f)^*_{D_0}$ relative to a domain $D_0$ in the 
following setup: $D\subset D_0$ are domains in $\CC$, $K=D_0\setminus D$ is a closed polar 
subset of $D$, and $z\in K$. We showed that a necessary and sufficient conditions for
$\{z\}\times\CC\cap (\Gamma_f)^*_{D_0\times \CC}=\emptyset$ is that $z$ be a regular boundary 
point for the Dirichlet problem on $D_M=\{\zeta\in D_0:\, |f(\zeta)|<M\}$.

In the present paper we continue our study of  $(\Gamma_f)^*_{D_0\times\CC}$. 
Our main results in that direction are Theorem \ref{thm:5.8} and Theorem \ref{thm:5.9} in 
Section 5, stating that if $(\{z\}\times\CC)\cap (\Gamma_f)^*_{D_0\times\CC}$ is not empty, 
then it consists of exactly one point. Thus a complete description is obtained
 of the pluripolar hulls of graphs of holomorphic functions that have a polar singularity set.  
As an important tool we introduce in this section the notion of {\it interior values} of 
holomorphic mappings.
These give rise to non-trivial points in the pluripolar hull of graphs of holomorphic mappings. 
In the one-dimensional case we show that the interior values of $f$ --- if they exist --- are unique 
and coincide with the value of a {\em distinguished homomorphism } as introduced by Gamelin and
Garnett, cf.~\cite{GG}. In \cite{EW1,EW2} we gave a sufficient condition for graphs of holomorphic functions to have a non-trivial pluripolar hull; Theorem \ref{thm:5.2} provides the natural generalization to pluripolar sets.

As a preparation we study in Section 2
pluriharmonic measure and extend work of Levenberg and Poletsky, \cite{LP},
as well as some results in \cite{EW2} on this topic. Noteworthy is Theorem \ref{thm:2.3}, which leads rapidly to the just mentioned Theorem \ref{thm:5.2}.  As one may expect, knowledge of pluriharmonic measure can be translated to pluripolar hulls. This is done in Section 3.
 
In Section 4 we prove a localization principle for pluriharmonic measure. This turns out to be
strong enough to explain qualitatively Siciak's \cite{Si} extension of our example in
\cite{EW1} of a holomorphic function $f\in A^{\infty}(\DD)$ with domain of
existence the unit disc $\DD$, which has $(\Gamma_f)^*$ extending over most of $\CC$.
We also show that the pluripolar hull of a connected $F_\sigma$-pluripolar 
set is connected; this may be of independent interest.

Throughout the paper $\BB(a,r)$ denotes the ball in $\CC^n$, centered at $a$ with radius $r$. 

The first named author thanks to Norm Levenberg for very helpful discussions, the second author is grateful to Tony O'Farrell for a useful conversation.


\section{Pluriharmonic measure} 
Let $\Omega$ be an open set in $\CC^n$ and let $E\subset\overline{\Omega}$
be any subset. The {\it pluriharmonic measure} of $E$ relative to $\Omega$
(or, {\it relative extremal function}) is defined as follows
(see e.g. \cite{K})
\begin{multline}\label{eq:1.1}
\omega(z,E,\Omega)=-\sup\{u(z):u\in\psh(\Omega), u\le0\text{ on } \Omega,\\
\text{ and }
\limsup_{\Omega\ni w\to\zeta} u(w)\le-1\text{ for }\zeta\in E\},
\quad z\in \Omega.
\end{multline}

Note that the function $-\omega$ need not be in PSH, but if $E$ is open then 
$-\omega\in \psh(\Omega)$.

Let $h: \Omega_1\to\Omega_2$ be holomorphic. A straightforward consequence of the definition is, see
\cite{LP},
\begin{equation}
\omega(z, h^{-1}(E),\Omega_1)\le \omega(h(z),E,\Omega_2), \quad z\in \Omega_1, E\subset \Omega_2.\label{eq:1.1a}
\end{equation}
\begin{prop}\label{prop:2.1} Let $\Omega'\Subset\Omega$ be open sets in $\CC^n$ and let
$V\Subset U\subset \Omega$ be open subsets.
Fix a $\zeta_0\in\overline{\Omega'}$.
Then there exists a neighborhood $W$ of $\zeta_0$ such that
\begin{equation}
\sup_{z\in W\cap\Omega'} \omega(z,V\cap\Omega',\Omega')\le\omega(\zeta_0,
U,\Omega).
\end{equation}
\end{prop}

\begin{proof}
There exists an $\eps>0$ such that $V_\eps\subset U$ and that
$\Omega'_\eps\subset\Omega$. Put $W=\BB(\zeta_0,\eps)$.
Then $V-w\subset U$ for any
$w\in\BB(0,\eps)$. So $\omega(\zeta_0+w,V\cap\Omega',\Omega')
=\omega(\zeta_0,V\cap\Omega'-w,\Omega'-w)\le \omega(\zeta_0,U,\Omega)$
for any $w\in\BB(0,\eps)$ such that $\zeta_0+w\in W\cap\Omega'$.
\end{proof}

Recall the following very useful result (see \cite{LP})
\begin{prop}\label{prop:2.2} Let $\Omega$ be an open set in $\CC^n$ and let $E\subset \Omega$
be any subset. Then
\begin{equation}
\omega(z,E,\Omega)=\inf\{\omega(z,U,\Omega):
U\subset\Omega\text{ is open and }E\subset U\}.
\end{equation}
\end{prop}

For a subset $E$ of $\CC^n$ and for a $\delta>0$ we put
$$E_\delta=\{z\in\CC^n:\dist(z,E)<\delta\}.$$

Combination of the above Propositions yields immediately the following.
\begin{theorem}\label{thm:2.3}
Let $\Omega'\Subset\Omega$ be open sets in $\CC^n$  and let
$E\subset\Omega'$ be a compact subset. Assume that a sequence
$\{z_n\}_{n=1}^\infty\subset\Omega'$ tends to
$\zeta_0\in \overline\Omega'$.
Then
\begin{equation}
\limsup_{n\to\infty} \omega(z_n, E_{\frac1n},\Omega')\le
\omega(\zeta_0,E,\Omega).\label{eqn:2.4}
\end{equation}
\end{theorem}
\begin{proof} 
Indeed, from Proposition~\ref{prop:2.1} we have 
\begin{equation}
\limsup_{n\to\infty} \omega(z_n, E_{\frac1n},\Omega')\le
\omega(\zeta_0,E_{\frac1m},\Omega)
\end{equation}
for any fixed $m$. Now apply Proposition~\ref{prop:2.2} to obtain \eqref{eqn:2.4}.
\end{proof}

Recall the following result (see e.g. \cite{K}, Corollary 4.5.11).
\begin{theorem}\label{thm:2.4}
Let $D$ be a hyperconvex domain in $\CC^n$ and let $K\subset D$ be a compact set.
Then $\omega(\cdot,K,D)$ is upper semi-continuous.
\end{theorem}

As a corollary of Theorem~\ref{thm:2.3} we next present a variant
of Theorem~\ref{thm:2.4} that gives a little less than upper semi-continuity,
 but is valid for arbitrary open sets in $\CC^n$.

\begin{corollary}\label{main-cor}
Let $\Omega'\Subset\Omega$ be open sets in $\CC^n$  and let
$K\subset\Omega'$ be a compact subset.
Then for any $\zeta_0\in \overline{\Omega'}$ we have
\begin{equation}
\limsup_{z\to\zeta_0} \omega(z, K,\Omega')\le \omega(\zeta_0,K,\Omega).
\end{equation}
\end{corollary}

\begin{proof} Note that $\omega(z,K,\Omega')\le
\omega(z,K_{\frac1n},\Omega')$ for any $n\in\NN$.
\end{proof}

Using similar methods we give an alternative proof of a result
of N.~Levenberg and E.~Poletsky.

\begin{corollary}[Levenberg-Poletsky \cite{LP}]
Let $\Omega'\Subset\Omega$ be open sets in $\CC^n$  and let
$E\subset\Omega'$ be a compact subset. Assume that
$V\subset \Omega\setminus E$ is an open set and that $\zeta_0\in V\cap\Omega'$.
Put $K=\overline{(\partial V)\cap\Omega'}$.
Then there exists a $\zeta\in K$ such that
\begin{equation}
\omega(\zeta_0, E,\Omega')\le \omega(\zeta,E,\Omega).
\end{equation}
\end{corollary}

\begin{proof}
Fix $n$ so large that $E_{\frac1n}\subset\Omega'$. We claim that there exists
a sequence $\{z_n\}\subset\partial V\cap\Omega'$ with
$\omega(z_n,E_{\frac1n},\Omega')>\omega(\zeta_0,E,\Omega')-\frac1n$.

Indeed, assume that for every $z\in\partial V\cap \Omega'$ we have
\begin{equation}
\omega(z,E_{\frac1n},\Omega')\le\omega(\zeta_0,E,\Omega').
\end{equation}
 Because $E$ is open, the function $-\omega(\cdot , E_{\frac1n},\Omega')$ is plurisubharmonic.
Therefore
\begin{equation*}
v(z)=\begin{cases}
-\omega(z,E_{\frac1n},\Omega'),&\text{ for }z\in\Omega'\setminus V,\\
\max\{-\omega(z,E_{\frac1n},\Omega'), -\omega(\zeta_0,E,\Omega')+\frac1n\},&
\text{ for }z\in V\cap\Omega'.
\end{cases}
\end{equation*}
is in $\psh(\Omega')$. 
We have $v\le0$ on $\Omega'$, $v\le-1$ on $E$. Hence,
\begin{equation*}
-\omega(\zeta_0,E,\Omega')\ge v(\zeta_0)\ge -\omega(\zeta_0,E,\Omega')+\frac1n,
\end{equation*}
a contradiction.

The conclusion is that there exists a subsequence $\{z_{n_k}\}$ converging
 to $\zeta\in\overline{\partial V\cap \Omega'}$
such that 
\begin{equation}
\omega(\zeta,E,\Omega)\ge\limsup_{k\to\infty}\omega(z_{n_k},
E_{\frac1{n_k}},\Omega')\ge\omega(\zeta_0,E,\Omega').
\end{equation}
\end{proof}

The next theorem is very important in our theory. It extends Theorem~3.7 in
\cite{EW2}. 
\begin{theorem}\label{thm:2.7} Let $D$ be a bounded open set in $\CC$
and let $\Delta\Subset D$ be a closed disc. Assume that
$K\subset\partial D$ is a compact polar set. Then for any $z_0\in K$ we have
$$
\limsup_{\zeta\to z_0}\omega(\zeta,\Delta,D)=
\inf\{\omega(z_0,\Delta,D\cup U): K\subset U\text{ open}\}.
$$
In particular, if $z_0\in K$ is a regular  boundary point of $D$
then
\begin{equation*}
\inf\{\omega(z_0,\Delta,D\cup U): K\subset U\text{ open}\}=0.
\end{equation*}
\end{theorem}

\begin{proof}
Observe that $\le$ is evident. For the inequality $\ge$ we take an open neighborhood $U$ of
$K$ and note that
for every $0<\eps\le 1$ the set
\begin{equation*}
F^\eps_U=\{ z\in \overline{D\cup U}:\ \limsup_{\zeta\to z}
\omega(\zeta,\Delta,D\cup U)\ge\eps \}
\end{equation*}
is a compact connected subset of $\overline{D\cup U}$ that contains $\Delta$.
Moreover, if $U_1\subset U_2$ then $F^\eps_{U_1}\subset F^\eps_{U_2}$.
We set $F^\eps=\cap_U F^\eps_{U}$. Then $F^\eps$ is a compact connected
subset of $\overline D$.

Now let
\begin{equation}\label{j1}
\limsup_{D\ni \zeta\to z_0}\omega(\zeta,\Delta,D)=\alpha\ge0.
\end{equation}
As $F^\eps\cap \partial D$ is a subset of the union of $K$
and the irregular boundary points of $D$, the set $F^\eps\cap \partial D$ is
 thin at $z_0$ and therefore totally disconnected.

To reach a contradiction, suppose that
\[
\inf\{\omega(z_0,\Delta,D\cup U): K\subset U\text{ open}\}=\eps>\alpha,
\]
then $z_0\in F^\eps$.

For any decreasing sequence $\{U_i\}$ of neighborhoods of $K$
with $K=\cap U_i$, the functions
$\omega(z,\Delta,D\cup U_i)$ form a decreasing sequence of harmonic functions on 
$D\setminus \Delta$,
 and hence converge uniformly on compact sets in $D\setminus \Delta$.
The limit function is $\omega(z,\Delta,D)$ and hence $\omega(z,\Delta,D)\ge\eps$ 
on $F^\eps\cap D$.
In view of (\ref{j1}) we infer that there is a neighborhood $V$ of $z_0$
such that $F^\eps\cap V\subset \partial D$. Thus $z_0$ is not in the component of $\Delta$,
 which is a
contradiction.
\end{proof}

\section{Properties of pluripolar hulls}
We commence by recalling two important definitions.
Let $\Omega$ be an open set in $\CC^n$ and let $E\subset\overline{\Omega}$
be a pluripolar subset. The {\it pluripolar hull} of $E$ in $\Omega$
is defined as
\begin{multline}
E^\ast_\Omega=\{z\in\Omega: \text{ for all }u\in\psh(\Omega): u|_E=-\infty\implies
u(z)=-\infty\}.
\end{multline}
For a pluripolar set $E$ in an open set $\Omega$ in $\CC^n$,
Levenberg and Poletsky define
the {\it negative pluripolar hull of $E$\/} in $\Omega$ as
\begin{equation*}
E_\Omega^{-}:=\{z\in \Omega: \text{ for all } u\in\psh(\Omega),
u\le0:u|_E=-\infty \implies u(z)=-\infty\}.
\end{equation*}
We extend the above definition to arbitrary pluripolar sets $E\subset\CC^n$ as
follows
\[
E^{-}_{\Omega}=\cap_{U\supset E\text{ open}} E^{-}_{\Omega\cup U}.
\]

We will use the following two important results from \cite{LP}.
\begin{theorem}\label{thm:3.1}
Let $\Omega$ be an open set in $\CC^n$ and let $E$ be
a pluripolar set in $\Omega$. Then
\begin{equation*}
E_\Omega^-=\{z\in \Omega: \omega(z,E,\Omega)>0\}.
\end{equation*}
\end{theorem}

\begin{theorem}\label{thm:3.2}
Let $\Omega$ be a pseudoconvex domain and let $E\subset\Omega$ be pluripolar.
Suppose that $\Omega=\cup_j\Omega_j$, where 
$\Omega_j\subset\Omega_{j+1}$ form an increasing sequence of relatively compact pseudoconvex 
subdomains of $\Omega$. Then 
$$E^\ast_{\Omega}=\cup_{j=1}^\infty (E\cap\Omega_j)^{-}_{\Omega_j}.$$
\end{theorem}

From Theorem~\ref{thm:3.1} it follows that for a compact
pluripolar set $K$ its negative pluripolar hull $K^-_{\Omega}$
is of $G_\delta$-type. And, therefore, if $\Omega$ is pseudoconvex
then $K^\ast_{\Omega}$ is of type $G_{\delta\sigma}$. Hence it is
a Borel set.

The following theorem is a high-dimensional version of Theorem~\ref{thm:2.7}

\begin{theorem}\label{thm:3.3}
Let $\Omega$ be an open set in $\CC^n$
and let $E\subset\Omega$ be any subset. Assume that
$F\subset\CC^n$ is a pluripolar set. Then
\[
\omega(z,E,\Omega)=
\inf\{\omega(z,E,\Omega\cup U): F\subset U\text{ open}\},\quad z\in
\Omega\setminus F^{-}_{\Omega}.
\]
\end{theorem}

\begin{proof}
Note that the inequality ''$\le$'' is trivial.

Fix a point $z_0\in\Omega\setminus F^{-}_{\Omega}$.
There exists a neighborhood $U$ of $F$ and a negative
plurisubharmonic function $h$ on $\Omega\cup U$ such that
$h=-\infty$ on $F$ and $h(z_0)\ne-\infty$.

Fix an $\eps>0$ and put
$U_\eps=\{z\in U: h(z)<-\frac{1}{\eps}\}$.
Note that $U_\eps\subset U$ is an open neighborhood of
$F$.
Let $v$ be a negative plurisubharmonic function on $\Omega$
such that $v\le-1$ on $E$. Consider the plurisubharmonic function
\[
v_\eps(z)=\begin{cases} \max\{v(z)+\eps h(z),-1\},&
z\in\Omega\setminus U_\eps,\\
-1,& z\in U_\eps.
\end{cases}
\]
Note that
\[
-v_\eps(z)\ge\omega(z,E,\Omega\cup U_\eps)\ge
\inf\{\omega(z,E,\Omega\cup U): F\subset U\text{ open}\},\quad z\in
\Omega\setminus F^{-}_{\Omega}.
\]
We let $\eps\to0$ and get the result.
\end{proof}

\begin{theorem}\label{thm:3.4} Let $\Omega$ be an open set in $\CC^n$ and let $\Omega'\Subset
\Omega$. Assume that $E\subset\Omega'$ is a compact pluripolar subset.
Then for any sequence $\{z_n\}_{n=1}^\infty\subset \Omega'$ such that
$\limsup_{n\to\infty}\omega(z_n,E,\Omega')>0$ and that $z_n\to w_0$ it follows
that $w_0\in E^{-}_{\Omega}$. Moreover, if $\Omega$ is pseudoconvex, then
$w_0\in E^\ast_\Omega$.
\end{theorem}

\begin{proof} Corollary~\ref{main-cor} gives $\omega(w_0, E, \Omega)>0$. Thus Theorem~\ref{thm:3.1} shows $w_0\in E^-_\Omega$. Next if $\Omega$ is pseudoconvex we apply Theorem \ref{thm:3.2} on a suitable exhaustion $\cup_j\Omega_j$ of $\Omega$ and find $w_0\in E^*_\Omega$. 
\end{proof}


\section{A localization principle}

The following localization principle is a main tool in our theory. 
Special cases of it
appear  in \cite{W2} and \cite{EW2}.
\begin{theorem}[A localization principle]\label{thm:4.1}
Let $\Omega\subset\CC^n$ be an open set
and let $E$ be an $F_\sigma$-pluripolar subset of $\Omega$.
Then for any open set $\Omega'\Subset \Omega$ and any open set
$U$ such that $\partial U\cap E^\ast_\Omega=\varnothing$ we have
\begin{equation}
\omega(z,E\cap U\cap \Omega',\Omega')=
\omega(z,E\cap U\cap \Omega',U\cap \Omega'),
\quad z\in U\cap \Omega'.\label{local}
\end{equation}
\end{theorem}

The proof will be  based on two lemmas. Their statement and proof are similar 
to work of Zeriahi (cf.~\cite{Z}, Lemme 2.1).
\begin{lemma}\label{lem:4.2}
Let $\Omega\subset\CC^n$ be an open set and let
$E\subset\Omega$ be a pluripolar subset. Assume that $F\subset E$,
$K\subset\Omega\setminus E^{\ast}_{\Omega}$ are compact subsets
and that $\Omega'\Subset\Omega$ is an open set.
Then for any number $N>0$ there exists a continuous negative plurisubharmonic
function $v$ on $\Omega'$ such that $v\le-N$ on $F\cap\Omega'$,
$v\ge-1$ on $K\cap\Omega'$.
\end{lemma}

\begin{proof}
Let $a\in K\subset\Omega\setminus E^\ast_{\Omega}$.
By the definition of $E^\ast_{\Omega}$ there exists a plurisubharmonic
function $u$ on $\Omega$ such that $u|_{E}=-\infty$ and $u(a)>-\infty$.
Put $M=\max_{z\in\Omega'\cup K}\{u(z)-u(a),0\}$.
Then the function
\begin{equation*}
v(z)=\frac{1}{2M+1}(u(z)-u(a))-\frac12,\quad z\in\Omega,
\end{equation*}
is a plurisubharmonic function on $\Omega$ with $v|_{E}=-\infty$,
$v(a)=-\frac12$ and $v\le0$ on $\Omega'$.

By the main approximation theorem for plurisubharmonic function
(see \cite{K}), there exists a decreasing sequence $\{v_j\}$ of continuous
plurisubharmonic functions on $\Omega'$ which tends pointwise to $v$.

Let $N>0$ be fixed. Dini's lemma on monotone decreasing sequences of continuous functions  provides us with a number $j_a>1$ such that
$v_{j_a}\le-N$ on $F$ and $v_{j_a}\le0$ on $\Omega'$.
Since $v_{j_a}$ is continuous on $\Omega'$ and since
$v_{j_a}(a)\ge v(a)=-\frac12>-1$, we may find a neighborhood $U_a$ of $a$
such that $v_{j_a}\ge-1$ on $U_a$.

Using a standard compactness argument, we construct a continuous
plurisubharmonic function $\widetilde v=\max\{v_{j_{a1}},\dots,v_{j_{am}}\}$
on $\Omega'$ such that $v\le0$ on $\Omega'$, $v\le-N$ on $F\cap\Omega'$,
and $v\ge-1$ on $K\cap\Omega'$.
\end{proof}

An immediate corollary of Lemma~\ref{lem:4.2} is
\begin{lemma}\label{lem:4.3}
Let $\Omega\subset\CC^n$ be an open set and let
$E\subset\Omega$ be an $F_\sigma$-pluripolar subset. Assume that
$K\subset\Omega\setminus E^{\ast}_{\Omega}$ is a compact subset
and that $\Omega'\Subset\Omega$ is an open set.
Then there exists a negative plurisubharmonic function $v$ on $\Omega'$
such that $v=-\infty$ on $E\cap\Omega'$, $v\ge-1$ on $K\cap\Omega'$.
\end{lemma}

\begin{proof}[Proof of Theorem~ \ref{thm:4.1}]
Fix an open set $\Omega'\Subset \Omega$.
Since $U\cap \Omega'\subset \Omega'$, we have the inequality ''$\ge$'' in
\eqref{local}.

Let us show the inequality ''$\le$''. Note that
$K:=\partial U\cap\overline{\Omega'}$ is a compact subset of $\Omega$.
According to Lemma~\ref{lem:4.3} there exists a plurisubharmonic function
$v$ on $\Omega'$ such that:
\begin{itemize}
\item $v\le0$ on $\Omega'$
\item $v=-\infty$ on $E\cap\Omega'$;
\item $v\ge-1$ on $K\cap\Omega'$.
\end{itemize}
Let $h\in\psh(\Omega'\cap U)$ be such that
$h\le-1$ on $E\cap U\cap \Omega'$ and $h\le0$ on $U\cap \Omega'$.
Fix an $\eps>0$. We consider the function
\begin{equation*}
v_\eps(z):=\begin{cases}
\max\{h(z)-\eps,\eps v(z)\}, &\text{if }z\in\Omega'\cap U,\\
\eps v(z), &\text{if } z\in \Omega'\setminus U.
\end{cases}
\end{equation*}
Note that $v_\eps$ is a negative plurisubharmonic function on $\Omega'$
which satisfies $v_\eps\le-1$ on $E\cap\Omega'$. Hence,
\begin{equation*}
-\omega(z,E\cap U\cap\Omega',\Omega')\ge v_\eps(z)\ge h(z)-\eps,
\quad z\in\Omega'\cap U.
\end{equation*}
Let $\eps\to0$. Then
$-\omega(z,E\cap U\cap \Omega',\Omega')\ge h(z)$, $z\in\Omega'\cap U$.
Therefore,
\begin{equation*}
\omega(z,E\cap U\cap \Omega',\Omega')\le \omega(z,E\cap U\cap\Omega',U\cap\Omega').
\end{equation*}
\end{proof}

Our next Proposition is an easy consequence of Theorem \ref{thm:4.1}. 
We do not know if the condition that $E$ is an $F_\sigma$ may be omitted.

\begin{prop}\label{prop:4.4} Let $\Omega$ be a pseudoconvex open set in $\CC^n$ and let
$E\subset\Omega$ be an $F_\sigma$-pluripolar subset. Assume that
$E$ is connected. Then $E^\ast_{\Omega}$ is also connected.
\end{prop}

\begin{proof} Assume that $E^\ast_\Omega\subset U_1\cup U_2=U$,
where $U_1,U_2$ are open sets such that $U_1\cap U_2=\varnothing$.
Since $E$ is connected, $E\subset U_1$ or $E\subset U_2$. Assume that
$E\subset U_1$.

Let $\Omega'\Subset\Omega$ be an open set.
Then
\[
\omega(z,E\cap \Omega',\Omega')=
\omega(z,E\cap U\cap \Omega',\Omega')=
\omega(z,E\cap \Omega',U\cap \Omega'),
\quad z\in U\cap \Omega'.
\]
Hence, $\omega(z,E\cap \Omega',\Omega')=0$ for $z\in U_2\cap\Omega'$.
Therefore, $(E\cap\Omega')^{-}_{\Omega'}\cap U_2=\varnothing$ and
$E^\ast_{\Omega}\cap U_2=\varnothing$. Here we used Theorem \ref{thm:3.2}.
\end{proof}

\begin{remark} Using Poletsky's theory \cite{P1}, \cite{P2}
of holomorphic discs one can give another proof of Proposition~\ref{prop:4.4}
\cite{P}. 
\end{remark}

Note that if $f$ is a holomorphic function on the unit disc $\DD$,
then its graph $(\Gamma_f)^\ast_{\CC^2}$ is a connected set and, therefore,
$\pi((\Gamma_f)^\ast_{\CC^2})$ is also connected, where
$\pi:\CC^2\ni(z,w)\to z\in\CC$ is the projection to the first
coordinate. In particular, the set $\pi((\Gamma_f)^\ast_{\CC^2})$
is not thin at any point of itself. Here, we show that in some
cases it cannot contain boundary points. We obtain this as a
corollary of the following more general result.

\begin{theorem}\label{thm:4.5} Let $m,n\in\NN$ and let $E\subset\CC^n$ be
an $F_\sigma$-pluripolar subset. Assume that $F:\CC^n\to\CC^m$ is
a holomorphic mapping such that
\begin{itemize}
    \item $F(E)\subset\BB_m$;
    \item $F(E^\ast_{\CC^n})\subset\overline{\BB_m}$.
\end{itemize}
Then $F(E^\ast_{\CC^n})\subset\BB_m$.
\end{theorem}

\begin{proof} Assume that $z_0\in E^\ast_{\CC^n}$ is such that
$F(z_0)\in\partial\BB_m$. Fix an $\eps>0$ and $r\in(0,1)$. Put
$U_{\eps}=F^{-1}(\BB_{1+\eps})$. Note that $U_\eps$ is an open
neighborhood of $E^\ast_{\CC^n}$. By Theorem~\ref{thm:4.1} we have for any
$R>1$
\begin{multline}
\omega\big(z_0,E\cap F^{-1}(\BB_r)\cap\BB_R,\BB_R)=
\omega\big(z_0,E\cap F^{-1}(\BB_r)\cap\BB_R,\BB_R\cap U_\eps)\\
\le\omega\big(F(z_0),\BB_r,\BB_{1+\eps}\big).
\end{multline}
Since $\eps>0$ is arbitrary, we get
\[
\omega\big(z_0,E\cap F^{-1}(\BB_r)\cap\BB_R,\BB_R)=0.
\]
So, $z_0\not\in (E\cap F^{-1}(\BB_r)\cap\BB_R)^{-}_{\BB_R}$. Since
$R>1$ is arbitrary, it follows that $z_0\not\in F\Big(\big((E\cap
F^{-1}(\BB_r)\big)^\ast_{\CC^n}\Big)$. From \cite{E} we obtain that
$z_0\not\in F(E^\ast_{\CC^n})$.
\end{proof}

\begin{remark} The first condition in Theorem~\ref{thm:4.5}
(i.e. $F(E)\subset\BB_m$) is essential. Indeed, in \cite{LMP}
a function $f\in\OO(\DD)\cap C^\infty(\overline\DD)$ is constructed such that
the graph
\[
\Gamma_f=\{(z,f(z)): z\in\overline\DD\}
\]
is complete pluripolar in $\CC^2$. Hence, $\pi((\Gamma_f)^\ast_{\CC^2})=\overline\DD$,
where $\pi$ is the projection.
\end{remark}

\begin{corollary}\label{cor:4.6} Let $f\in\OO(\DD)$ be a holomorphic function such that
$(\Gamma_f)^\ast_{\CC^2}\subset \overline{\DD}_\rho\times\CC$,
where $\rho\ge1$. Then $(\Gamma_f)^\ast_{\CC^2}\subset
{\DD_\rho}\times\CC$.
\end{corollary}

In \cite{EW1}, the authors constructed an example of a smooth
holomorphic function $f$ on the unit disc such that
$(\Gamma_f)^\ast_{\CC^2}\setminus\Gamma_f\ne\varnothing$. From
Proposition~\ref{prop:4.4} (see the discussion after the Proposition) and 
Corollary~\ref{cor:4.6} we see that the set
$(\Gamma_f)^\ast_{\CC^2}\setminus\Gamma_f$ is actually quite big. See also Siciak
\cite{Si}.

\begin{corollary}\label{cor:4.7}
Let $f\in\OO(\DD)$ be a holomorphic function. Assume that $r_n\searrow1$
is a sequence of radius such that
$(\Gamma_f)^\ast_{\CC^2}\cap (\partial\DD_{r_n}\times\CC)=\varnothing$. Then
$(\Gamma_f)^\ast_{\CC^2}\subset\DD\times\CC$. Moreover, if $f$ is bounded then 
$\Gamma_f$ is complete pluripolar.
\end{corollary}

\begin{proof} The first part follows from Corollary~\ref{cor:4.6}. So, assume that
$f$ is bounded. Fix a closed disc $S\subset\DD$, denote the graph of $f$ over $S$ by $\Gamma_S$,
 and put
$R_0=\max\{1,\sup_{\DD}|f|\}$.  Then for any $R>R_0$ we have
$(\Gamma_S)^-_{\DD_R\times\DD_R}\cap (\DD\times\DD_R)=
(\Gamma_S)^{-}_{\DD\times\DD_R}=\Gamma_f$. Hence, $(\Gamma_f)^\ast_{\CC^2}=
(\Gamma_S)^\ast_{\CC^2}=\cup_{R>R_0}(\Gamma_S)^-_{\DD_R\times\DD_R}=\Gamma_f$.
\end{proof}

As a simple corollary of the localization principle we have the following
\begin{corollary}
Let $\Omega$ be a pseudoconvex domain in $\CC^n$ and let $E\subset\Omega$
be an $F_\sigma$-pluripolar set such that $E^\ast_\Omega\Subset\Omega$. Then for any
open set $\Omega'\Subset\Omega$ such that $E^\ast_{\Omega}\subset\Omega'$ we have
$E^-_{\Omega'}=E^\ast_\Omega$.
\end{corollary}

\begin{proof} Let $\Omega''$ be a pseudoconvex domain such that 
$\Omega'\Subset\Omega''\Subset\Omega$. From the localization principle we have
$\omega(z,E,\Omega'')=\omega(z,E,\Omega')$ for $z\in\Omega'$. Hence, 
$E^-_{\Omega''}=E^-_{\Omega'}$. Since $\Omega''$ is arbitrary,
$E^\ast_\Omega=E^-_{\Omega'}$.
\end{proof}

\section{The set of interior values}
In the study of boundary behavior of a holomorphic function
the properties of its cluster set are very important
(see e.g. \cite{N}).
In connection with the pluripolar hull a certain subset of the
cluster set is very useful.

\begin{definition}
Let $\Omega\subset\CC^n$ be an open set and let $f:\Omega\to\CC^m$
be a bounded holomorphic mapping. Assume that $z_0$ is a boundary
point of $\Omega$.
An {\em interior value} of $f$ at $z_0$
is a limit point of a sequence $f(z_k)$, where $z_k\in \Omega$ tend to $z_0$
in such a way that for some closed non-empty
ball $B\subset\Omega$ and some positive
number $\alpha$ we have
\begin{equation}\label{eq:5.1}
\omega(z_k,B,\Omega)\ge \alpha\quad \text{ for any } k\ge1.
\end{equation}

 We denote the set of interior values of $f$
at $z_0\in\partial\Omega$ by $L_{z_0}(f;\Omega)$. 

For an unbounded holomorphic mapping $f$ defined on an open set $\Omega$
we put $L_{z_0}(f;\Omega)=\cup_{R>0} L_{z_0}(f;\Omega_R)$, where
\[
\Omega_R=\{z\in \Omega: \|f(z)\|<R\}
\]
and $\|\cdot\|$ denotes the Euclidean norm in $\CC^m$.
\end{definition}
In case 
for some (and, therefore, for any) closed non-empty
ball $B\subset\Omega$ we have
$\lim_{z\to z_0}\omega(z,B,\Omega)=0$, we put
$L_{z_0}(f;\Omega)=\varnothing$.
This happens if and only if $z_0$ is a regular boundary point of $\Omega$ for the 
Dirichlet problem.

The following little lemma shows that in $\CC$  interior
value is a ''local property''.

\begin{lemma} \label{lem:5.2a}Let $\Omega$ be an open set in $\CC$ and let
$z_0\in\partial \Omega$. Assume that $f:\Omega\to\CC^m$ is a holomorphic
mapping. Then there exists an $r>0$ such that
\begin{equation}
L_{z_0}(f;\Omega)=L_{z_0}(f;\Omega\cap \DD(z_0,r)),
\end{equation}
where $\DD(z_0,r)=\{z\in\CC: |z-z_0|<r\}$.
\end{lemma}

\begin{proof} This follows from Bouligand's lemma (see \cite{R}).
\end{proof}

\begin{theorem}\label{thm:5.2} Let $\Omega$ be an open set in $\CC^n$ and let $\Omega'\Subset
\Omega$. Assume that $E\subset\Omega'$ is a compact pluripolar subset.
Then for any sequence $\{z_n\}_{n=1}^\infty\subset \Omega'$ with $z_n\to w_0$
and such that
$\limsup_{n\to\infty}\omega(z_n,E,\Omega')>0$, it follows
that $w_0\in E^{-}_{\Omega}$. Moreover, if $\Omega$ is pseudoconvex, then
$w_0\in E^\ast_\Omega$.
\end{theorem}

\begin{proof} This is a direct consequence of Corollary~\ref{main-cor}.
\end{proof}

\begin{corollary}\label{intval1}
Let $\Omega'\Subset\Omega$ be open sets in $\CC^n$ and let
$f:\Omega'\to\CC^m$ be a holomorphic mapping. Assume that
$\zeta_0\in\partial\Omega'$. Then $\{\zeta_0\}\times L_{\zeta_0}
(f;\Omega')\subset (\Gamma_f)^\ast_{\Omega\times{\CC^m}}$.

In particular, if $L_{\zeta_0}(f;\Omega')$ is non-pluripolar, then
\[
(\Gamma_f)^\ast_{\Omega\times\CC^m}\cap \{\zeta_0\}\times\CC^m=
\{\zeta_0\}\times\CC^m.
\]
\end{corollary}

For $n=1$ we have a little bit stronger result.
\begin{corollary}\label{intval2}
Let $\Omega'\subset\Omega$ be open sets in $\CC$ and let
$f:\Omega'\to\CC$ be a holomorphic function. Assume that
$\zeta_0\in\partial\Omega'\cap\Omega$.
Then $\{\zeta_0\}\times L_{\zeta_0}(f;\Omega')\subset
(\Gamma_f)^\ast_{\Omega\times\CC}$.
\end{corollary}
\begin{proof}[Proof of both corollaries] Fix a ball $B\subset\Omega'$.
The inequality \eqref{eq:1.1a} with the map $h: z\mapsto(z,f(z))$ and the estimate \eqref{eq:5.1}
provide us with points  $w_k=(z_k,h(z_k))\in\Gamma_f\cap \Omega'\times \CC^m$ converging to $(\zeta_0,\eta_0)\in \{\zeta_0\}\times L_{\zeta_0}(f,\Omega')$ such that
$\limsup_{k\to\infty}\omega(w_k,h(B),\Omega'\times\BB(\eta_0,R))>0$. 
Now Theorem \ref{thm:5.2} 
applies. Use Lemma \ref{lem:5.2a} for Corollary \ref{intval1}.
\end{proof}

Let $D$ be a domain in $\CC$ and let $f\in\OO(D)$. Assume that $z_0\in\partial D$.
We want to show that $\#L_{z_0}(f;D)\le1$ and, therefore, the set
$ L_{z_0}(f;D)$ is always polar.
The crucial ingredient is work of Gamelin and Garnett \cite{GG}, which extends earlier work of Zalcman \cite{Za}. We recall it here for a small part.
 Consider $H^\infty(D)$, the algebra of bounded holomorphic functions on $D$. A {\em distinguished homomorphism} at $z_0$ is a homomorphism
above $z_0$ that belongs to the same Gleason part of the maximal ideal space
$\calm$  of $H^\infty(D)$ as the point evaluations at points of $D$.
Distinguished homomorphisms need not exist, but it is shown in \cite{GG} that there can at most be one distinguished homomorphism above $z_0$.

\begin{lemma}\label{uniqueness}
Let $D$ be a domain in $\CC$ and let $f\in\OO(D)$.
Assume that $z_0\in\partial D$. Then $\#L_{z_0}(f;D)\le1$.
\end{lemma}

We sketch one proof here and  give another one later.

\begin{proof}
[Sketch] We may assume $f$ is bounded on $D$. Let $B(a,r)$ be a compact ball in $\Omega$ 
and $\{z_n\}$ a sequence in $D$ tending to $z_0$ such that
\begin{equation}\label{eq:5.2}
\omega(z_k, B, D)\ge \alpha>0.
\end{equation}
Let $\phi_{z_n}$ be the associated point evaluations. The pseudo-hyperbolic distance is
$$d(\phi_a,\phi_{z_n})=\sup_{\{f\in H^\infty(D) :\ |f|<1, f(a)=0\}}\{|f(z_n)|\}.$$
By the two constant theorem and \eqref{eq:5.2}, $d(\phi_a,\phi_{z_n})\le c<1$.
Therefore any limit point $\mu$ of $\phi_{z_n}$ in $\calm$ has $d(\phi_a,\mu)\le c$ so 
such a $\mu$ is a distinguished homomorphism and must be unique.
Hence also $\lim f(z_n)=\mu(f)$ exists independently of the sequence $z_n$ with 
\eqref{eq:5.2}.
\end{proof}
\begin{remark}
It is well possible that a regular boundary point admits a distinguished homomorphism. Existence of distinguished homomorphisms can be characterized in terms of analytic capacity (Melnikov type condition), cf.~\cite{GG}, while regularity is characterized in terms of logarithmic capacity (Wiener's criterion), cf.~\cite{R}.
\end{remark}

The other proof will be based on the connection between distinguished
 homomorphisms and interpolating sequences.
A sequence $\{z_n\}_{n=1}^\infty\subset D$ is called an {\em interpolating sequence }
for $H^\infty(D)$ if for every bounded sequence $\{s_n\}_{n=1}^\infty\in\ell^\infty$,
there is $f\in H^\infty(D)$ such that $f(z_n)=s_n$ for any $n\ge1$.

Let us show the following variation of the well-known result
related to the Green function $g_D$ of a domain $D$ 
(see e.g. \cite{R}, Corollary 4.5.5).
\begin{prop}\label{interp}
Let $D$ be a domain in $\CC$ and let $\{z_n\}_{n=1}^\infty\subset D$ be an interpolating sequence. Then $\lim_{n\to\infty} g_D(z_n,z_1)=0$.
\end{prop}

\begin{proof} There exists a bounded holomorphic function $f$ on $D$ such that
$f(z_1)=1$ and $f(z_n)=0$ for any $n\ge 2$. Assume that $\|f\|=M$. Then
$g_D(z; \{z_n\}_{n=2})\ge \log |f(z)|-\log M$ and, therefore,
\[
\sum_{n=2}^\infty g_D(z_1;z_n)=g_D(z_1; \{z_n\}_{n=2}^\infty)\ge -\log M.
\]
Hence, $g_D(z_n,z_1)=g_D(z_1,z_n)\to0$ when $n\to\infty$.
\end{proof}

\begin{prop}\label{main-prop}
Let $D$ be a domain in $\CC$ and let $f\in\OO(D)$.
Assume that $z_0\in\partial D$ is an irregular point. Then there exists
$w_0\in\CC$ such that for any sequence $\{z_n\}\subset D$ with $z_n\to z_0$
there exists a subsequence $\{z_{n_k}\}$ such that $f(z_{n})\to w_0$ or
$g_D(z_{n_k},z_1)\to 0$.
\end{prop}

\begin{proof} From Theorem 4.5 in \cite{GG} we infer that there exists
$w_0\in\CC$ (the value of the distinguished homomorphism at $f$) such that
$f(z_n)\to w_0$ or there exists an interpolating subsequence of $\{z_n\}$.
Now the result follows from Proposition~\ref{interp}.
\end{proof}

\begin{proof}[second proof of Lemma \ref{uniqueness}]
 In case $f$ is bounded
Proposition~\ref{main-prop} applies. The general case follows from the
definitions.
\end{proof}

\begin{theorem}\label{thm:5.8}
Let $D$ be an open set in $\CC$ and let $A\subset D$
be a closed polar set. Assume that $f\in\OO(D\setminus A)$ and
that $z_0\in A$.
Then
\[
(\Gamma_f)_{D\times\CC}^\ast\cap\{z_0\}\times\CC=\{z_0\}\times L_{z_0}(f;D).
\]
And, therefore, $\#\Big((\Gamma_f)^\ast_{D\times\CC}\cap \{z_0\}\times\CC\Big)\le1$.
\end{theorem}

For the proof, first let us show the following refinement of the main
result of \cite{EW2}.
\begin{theorem}\label{thm:5.9}
Let $D$ be an open set in $\CC$ and let $A$ be a closed polar subset of $D$.
Suppose that $f\in\OO(D\setminus A)$ and that $z_0\in A$.
Assume that $U\subset\CC$ is an open set.
Then the following conditions are equivalent:
\begin{enumerate}
\item\label{eq3:1} $(\{z_0\}\times\CC)\cap (\Gamma_f\cap (D\times U))^\ast_{D\times U}=\varnothing$;
\item\label{eq3:2} there exists a sequence of open sets $V_1\subset  V_2\subset\dots\Subset U$
such that $\cup_j V_j=U$ and the set $\{z\in D\setminus A: f(z)\in U \setminus{V_j}\}$
is not thin at $z_0$ for any $j\ge1$.
\item\label{eq3:3} for any open set $V\Subset U$ the set
$\{z\in D\setminus A: f(z)\in U\setminus {V}\}$ is not thin at $z_0$.
\end{enumerate}
Moreover, if the set $\{z\in D\setminus A: f(z)\not\in V\}$
is thin at $z_0$ for some open set $V\Subset U$, then there exists a $w_0\in\overline{V}$,
such that $(z_0,w_0)\in(\Gamma_f\cap D\times U)^\ast_{D\times U}$.
\end{theorem}

\begin{proof} 

$(1)\implies (3)$. Assume that there exists an open set $V\Subset U$ such that 
$\{z\in D\setminus A: f(z)\not\in V\}$ is thin at $z_0$. Then the set
 $\{z\in D\setminus A: f(z)\in V\}$ is not regular at $z_0$. Hence, there exist
an open set $G\Subset D$ such that $\partial G\cap A=\varnothing$, $z_0\in G$, and a sequence
$\{z_n\}_n$ in $G\setminus A$ tending to $z_0$ such that
 $\limsup_{n\to\infty}\omega(z_n,S,G\setminus A)>0$ for some closed disc
$S\subset G\setminus A$. There is a subsequence $\{z_{n_k}\}$ such that
$f(z_{n_k})$ converges to an interior value $w_0\in\overline{V}$ and, using Theorem \ref{thm:5.2}
$(z_0,w_0)\in (\Gamma_f)^\ast_{D\times U}$. We have also proved the last statement of the theorem.

$(3)\implies(2)$. Obvious.

$(2)\implies(1)$. 
Again $\Gamma_S$ will denote the graph of $f$ over a disc $S$ in $D$.
In view of Theorems \ref{thm:3.1} end \ref{thm:3.2}, it suffices to show that for $w\in V_j$
$\omega((z_0, w), \Gamma_S, G\times V_j)=0$ for any fixed, open set $G\Subset D$ such that $\partial G\cap A=\varnothing$ and some
closed disc $S\subset G\setminus A$. 
 To estimate $\omega((z,f(z)), \Gamma_S, G\times V_j)$, let $\eps>0$ and start with a small neighborhood $V$ of $A\cap G$, to be determined later.
 Put $\widetilde V=V\cup (D\setminus\overline{G})$. Let
\begin{equation*}
U=[\{z\in D\setminus A: f(z)\in V_{j+1}\}\cup\widetilde V]\times V_{j+1}\cup
\{z\in D\setminus A: f(z)\not\in\overline{V_j}\}\times (\CC\setminus\overline{V_j}).
\end{equation*}
Then $U$ is a neighborhood of $\Gamma_f$.
It was proved in \cite{EW2} that 
\begin{equation}
(\Gamma_f)^*_{D\times\CC}\subset \Gamma_f\cup A\times\CC.\label{eq:5.11a}
\end{equation} Therefore $\partial U\cap (\Gamma_S)^*_{G\times V_j}=\emptyset$. We may apply the
the localization principle,  Theorem \ref{thm:4.1} and find
\begin{multline}
\omega((z,w),\Gamma_S,G\times V_j)=\omega((z,w),\Gamma_S,G\times V_j\cap U)\\
=
\omega((z,w),\Gamma_S,[\{z\in D\setminus A: f(z)\in V_{j+1}\}\cup\widetilde V]\times V_j),\label{eq:5.12a}
\end{multline}
for $(z,w)\in U\cap G\times V_j)$. Now we apply \eqref{eq:1.1a} to the projection $(z,w)\mapsto z$
and find that the right-hand side of \eqref{eq:5.12a} is 
\begin{equation}
\le
\omega(z,S, \{z\in D\setminus A: f(z)\in V_{j+1}\}\cup\widetilde V).
\end{equation}
By Theorem \ref{thm:2.7} we can choose $V$ so small that $ \omega(z_0,S, \{z\in D\setminus A: f(z)\in V_{j+1}\}\cup\widetilde V)<\eps$.
Letting $\eps \to 0$, it follows that $\omega((z_0,w),\Gamma_S,G\times V_j)=0$.
\end{proof}

For the proof of the main result we need the following simple remark related to
the pluripolar hull.

\begin{lemma}\label{lemma2}
Let $D\subset\CC^n$ be a pseudoconvex set and let $A\subset D$
be a closed pluripolar subset. Assume that $E\subset D\setminus A$
is a pluripolar compact set. Then
$E^\ast_{D}\subset E^\ast_{D\setminus A}\cup A$.
\end{lemma}

\begin{proof} Let $D_1\Subset D_2\Subset\dots\Subset D$ be an exhaustion
of $D$ by hyperconvex domains. Then by Theorem \ref{thm:3.2}
we have $E^\ast_D=\cup_{j=1}^\infty (E\cap D_j)^-_{D_j}=
\cup_{j=1}^\infty E_{D_j}^-$. Now we apply
Lemma 3.1 from \cite{LP}, saying that $\omega(z, E, D_j)=\omega(z,E, D_j\setminus A)$,
  and infer that $E^-_{D_j}\cap D_j\setminus A=
E^-_{D_j\setminus A}$, and the lemma follows.
\end{proof}

\begin{proof}[Proof of Theorem~\ref{thm:5.8}]
Assume that $L_{z_0}(f;D)\subset\{w_0\}$.
Put $U:=\CC\setminus\{w_0\}$. Then, by the definition of interior value, for every relative compact subset $W\Subset U$ the set
$\{z\in D\setminus A:\ f(z)\in U\setminus W\}$ is not thin at $z_0$.
Hence by Theorem \ref{thm:5.9}
$\{z_0\}\times U\cap (\Gamma_S)^\ast_{D\times U}=\varnothing$.
But $(\Gamma_S)^\ast_{D\times\CC}\subset (\Gamma_S)^\ast_{D\times U}\cup
(D\times \{w_0\})$. Therefore,
$\{z_0\}\times\CC\cap (\Gamma_S)^\ast_{D\times\CC}\subset\{(z_0,w_0)\}$.
\end{proof}

\begin{remark} Let $A=\{a_n\}_{n=1}^\infty\subset\DD\setminus\{0\}$ be a sequence
such that $a_n\to0$ and let $\{c_n\}_{n=1}^\infty\subset\CC\setminus\{0\}$.
Put 
\begin{equation}
f(z)=\sum_{n=1}^\infty\frac{c_n}{z-a_n}.\label{eq:5.11}
\end{equation}
 Suppose that
$\sum_{n=1}^\infty|c_n|<+\infty$ and that $\sum_{n=1}^\infty\left|\frac{c_n}{a_n}\right|$ converges.
Then $f\in\OO(\CC\setminus (A\cup\{0\}))$ and $f(0)$ is well-defined.

In \cite{EW2} the authors gave sufficient conditions on $\{a_n\}$ and
$\{c_n\}$ such that $(\Gamma_f)^\ast_{\CC^2}=\Gamma_f\cup\{(0,f(0)\}$.

Theorem~\ref{thm:5.8} gives that 
$\#\Big((\Gamma_f)^\ast_{\CC^2}\setminus \Gamma_f\Big)\le1$. 
In case $(\Gamma_f)^\ast_{\CC^2}=\Gamma_f\cup\{(0,w_0)\}$ it seems likely that $w_0=f(0)$, as defined by the series. Under mild convergence conditions this is easily proved.
\end{remark}

\begin{example} Suppose that the series \eqref{eq:5.11} has the property that
 $(\Gamma_f)^\ast_{\CC^2}$ contains a point $(0,w_0)$ and suppose that for every $M$ either 
the series
\begin{equation}\sum_{n=1}^\infty\left|\frac{c_n}{z-a_n}\right|\label{eq:5.12}
\end{equation}
is bounded on $\{z:\ |f(z)|<M\}$,
or the function 
\begin{equation}g(z)=\sum_{n=1}^\infty \frac{c_n}{a_n(z-a_n)}\label{eq:5.13}
\end{equation}
is in $H^\infty(D_M)$. Then $f(0)=w_0$.

From \cite{GG} we know that the distinguished homomorphism $\phi_0$ at 0 can be represented by a positive measure $\mu_0$ on $D_m$. Note that $\phi_0(c_n/(z-a_n))=-c_n/a_n$, because $c_n/(z-a_n)$ is holomorphic in a neighborhood of 0. If \eqref{eq:5.12} is satisfied, 
then by the dominated convergence theorem
$$f(0)=\lim_{N\to\infty}\sum_{n=1}^N
\frac{-c_n}{a_n}=\lim_{n\to N}\sum_{n=1}^N\int_{D_M} \left(\frac{c_n}{z-a_n}\right)d\mu_0
=\int_{D_M}f\, d\mu_0=w_0.$$ 
In case of \eqref{eq:5.13} we observe that
$$f(z)-f(0)=\sum_{n=1}^\infty \frac{zc_n}{a_n(z-a_n)}=zg(z).$$
Hence $\phi_0(f(z)-f(0))=0$, or $f(0)=w_0$.
\end{example}

%

\bibliographystyle{amsplain}


\end{document}